\documentclass[3p,times]{elsarticle}
\usepackage{amssymb}
\usepackage{amsmath}
\usepackage{mathrsfs}
\numberwithin{equation}{section}
\newtheorem{thm}{Theorem}[section]

\newtheorem{prop}[thm]{Proposition}
\newdefinition{defi}[thm]{Definition}
\newproof{pf}{Proof}
\usepackage[figuresright]{rotating}

\journal{Statistics and Probability Letters (Published in 2016)}

\begin{document}

\begin{frontmatter}

\title{On quasi-ergodic distribution for one-dimensional diffusions}

\author{Guoman He\corref{cor1}}
\ead{hgm0164@163.com}
\author{Hanjun Zhang}
\ead{hjz001@xtu.edu.cn}
\cortext[cor1]{Corresponding author}
\address{School of Mathematics and Computational Science, Xiangtan University, Hunan 411105, PR China}

\begin{abstract}
In this paper, we study quasi-ergodicity for one-dimensional diffusion $X$ killed at 0, when 0 is an exit boundary and $+\infty$ is an entrance boundary. Using the spectral theory tool, we show that if the killed semigroup is intrinsically ultracontractive, then there exists a unique quasi-ergodic distribution for $X$. An example is given to illustrate the result. Moreover, the ultracontractivity of the killed semigroup is also studied.
\end{abstract}


\end{frontmatter}

\section{Introduction}
\label{sect1}
One of the fundamental problems for a killed Markov process conditioned on survival is to study its long-term asymptotic behavior. In order to understand the behavior of the process before absorption, a relevant object to look at is a so-called quasi-ergodic distribution, which may be defined as the limiting distribution of the expected empirical distributions of the process conditioned on absorption not occurred. This paper is a continuation of studying quasi-ergodic distribution for killed Markov processes, which focus on one-dimensional diffusions with 0 as an exit boundary and $+\infty$ as an entrance boundary.
\par
Quasi-ergodic distribution is sometimes called mean-ratio quasi-stationary distribution in some literature. To the best of our knowledge, the definition of quasi-ergodic theorem put first forward by Breyer and Roberts \cite{BR99}, who proved the existence and uniqueness of quasi-ergodic distribution for Markov processes with general state spaces under the condition that the process is a positive Harris  $\lambda$-recurrent process with $\lambda\leq0$. For absorbing Markov processes under the assumption on $\lambda$-positivity, Chen et~al. \cite{CLJ12} studied some problems related to quasi-ergodic distribution, attempting to interpret quasi-ergodic distribution from different perspectives. Soon afterwards, Chen and Jian \cite{CJ13} proved the existence and uniqueness of both quasi-stationary distribution and quasi-ergodic distribution for killed Brownian motion by using an eigenfunction expansion for the transition density. Their results show that quasi-ergodic distribution is quite different from quasi-stationary distribution, but coincide with a certain double limiting distribution, reflecting a kind of phase transition. The exponential convergence rate to quasi-stationarity and quasi-ergodicity for ultracontractive Markov processes was also studied by Chen and Jian \cite{CJ14}.
Recently, for general Markov processes, which is a standard Markov process admitting a dual with respect to a finite measure and a strictly positive bounded continuous transition density(in fact, this condition is equivalent to saying that the semigroup is ultracontractive), Zhang et~al. \cite{ZLS14} proved the existence and uniqueness of both quasi-stationary distribution and quasi-ergodic distribution. We point out that our case cannot be contained as one of their cases because the reference measure is infinite in our case.
\par
This paper is related to \cite{LJ12}. Littin \cite{LJ12} proved the existence of a unique quasi-stationary distribution and of the Yaglom limit for the one-dimensional diffusion $X$ killed at 0, when 0 is an exit boundary and $+\infty$ is an entrance boundary. Using the spectral theory tool, we will show that if the killed semigroup is intrinsically ultracontractive, then there exists a unique quasi-ergodic distribution for the process $X$.
\par
The remainder of this paper is organized as follows. In Section \ref{sect2} we present some preliminaries that will be needed in the sequel. In Section \ref{sect3}, we will study some problems related to quasi-ergodic distribution and prove the existence and uniqueness of quasi-ergodic distribution. We conclude in Section \ref{sect4} with an example.

\section{Preliminaries}
\label{sect2}
We consider the generator $Lu:={1\over2}\partial_{xx}u-\alpha\partial_xu$. Denote by $X$ the diffusion on $(0,\infty)$ whose infinitesimal generator is $L$, or in other words the solution of the stochastic differential equation (SDE)
\begin{equation}
\label{2.1}
dX_t=dB_t-\alpha(X_t)dt,~~~~~~~~~~X_{0}=x>0,
\end{equation}
where $(B_t;t\geq0)$ is a standard one-dimensional Brownian motion and $\alpha\in C^1(0,\infty)$. In this paper, $\alpha$ is allowed to explode at the origin. There exists a pathwise unique solution to the SDE (\ref{2.1}) up to the explosion time $\tau$. \par
 Associated with $\alpha$, we consider the following two functions
\begin{equation}
\label{2.2}
\Lambda(x)=\int_{1}^{x}e^{Q(y)}dy~~~~~\mathrm{and}~~~~~\kappa(x)=\int_{1}^{x}e^{Q(y)}\left(\int_{1}^{y}e^{-Q(z)}dz\right)dy,
\end{equation}
where $Q(y)=\int_{1}^{y}2\alpha(x)dx$. Notice that $\Lambda$ is the scale function for $X$.
\par
The other important piece of information is the following measure, which is an infinite measure in this paper, defined on $(0,\infty)$:
\begin{equation}
\label{2.3}
\mu(dy):=e^{-Q(y)}dy.
\end{equation}
Notice that $\mu$ is the speed measure for $X$.
\par
Let $T_a:=\inf\{0\leq t<\tau:X_t=a\}$ be the hitting time of $a\in(0,\infty)$ for $X$. We denote by $T_{\infty}=\lim\limits_{n\rightarrow\infty}T_{n}$ and $T_{0}=\lim\limits_{n\rightarrow\infty}T_{1/n}$. Because $\alpha$ is regular in $(0,\infty)$, then $\tau=\min\{T_{0}, T_{\infty}\}$. Let $\mathbb{P}_x$ and $\mathbb{E}_x$ stand for the probability and the expectation, respectively, associated with $X$ when initiated from $x$. For any distribution $\nu$ on $(0,\infty)$, we define $\mathbb{P}_\nu(\cdot):=\int_{0}^{\infty}\mathbb{P}_x(\cdot)\nu(dx)$. We denote by $\mathcal{B}(0,\infty)$ the Borel $\sigma$-algebra on $(0,\infty)$, $\mathscr{P}(0,\infty)$ the set of all probability measures on $(0,\infty)$, ${\bf1}_{A}$ the indicator function of $A$ and $\langle f,g\rangle_\mu=\int_{0}^{\infty}f(u)g(u)\mu(du)$.
\par
For most of the results in this paper we will use the following hypothesis $(\mathrm{H})$, that is,
\begin{defi}
\label{defi1.1}
We say that hypothesis $(\mathrm{H})$ holds if the following explicit conditions on $\alpha$, all together, are satisfied:
\par $(\mathrm{H1})$ for all $x>0$, $\mathbb{P}_x(\tau=T_{0}<T_{\infty})=1$.
\par $(\mathrm{H2})$ for any $\varepsilon>0$, $\mu(0, \varepsilon)=\infty$.
\par $(\mathrm{H3})$ $S=\int_{1}^{\infty}e^{Q(y)}\left(\int_y^{\infty}e^{-Q(z)}dz\right)dy<\infty$.
\end{defi}
\par If (H1) holds, then it is equivalent to
$\Lambda(\infty)=\infty$ and $\kappa(0^{+})<\infty$
(see, e.g., \cite{IW89}, Chapter VI, Theorem 3.2). If (H1) and (H2) are satisfied, we say that $0$ is an exit boundary in the sense of Feller (see \cite[Chapter 15]{KT81}). If (H1) and (H3) are satisfied, then $+\infty$ is an entrance boundary in the sense of Feller ( see also \cite[Chapter 15]{KT81}).
\par
We know from \cite{LJ12} that $L$ is the generator of a strongly continuous symmetric semigroup of contractions on $\mathbb{L}^2(\mu)$ denoted by $(P_t)_{t\geq0}$. This semigroup is sub-Markovian, that is, if $0\leq f\leq1$, then $0\leq P_tf\leq1$ $\mu$-a.e.. Also from \cite{LJ12} we get that when absorption is sure, that is (H1) holds, the semigroup of $X$ killed at 0 can be given by $P_tf(x)=\mathbb{E}_x[f(X_t), T_{0}>t]$.
\par
In this paper, we need some results obtained by Littin \cite{LJ12}, which are summarized in the following proposition.
\begin{prop}
\label{prop 2.2}
$([10])$
Assume $(\mathrm{H})$ holds. Then we have
\par $(\mathrm{i})$ $-L$ has purely discrete spectrum. The eigenvalues $0<\lambda_1<\lambda_2<\cdots$ are simple, $\lim\limits_{n\rightarrow\infty}\lambda_n=+\infty$, and the eigenfunction $\eta_n$ associated to $\lambda_n$ has exactly $n$ roots belonging to $(0,\infty)$ and an orthonormal basis of $\mathbb{L}^2(\mu)$. In particular, $\eta_1$ can be chosen to be strictly positive.
\par
for $g\in \mathbb{L}^2(\mu)$,
\begin{equation}
\label{2.4}
P_tg=\sum_{i\geq1}e^{-\lambda_it}\langle\eta_i,g\rangle_\mu\eta_i~~~in~~~\mathbb{L}^2(\mu),
\end{equation}
then for $f,g\in \mathbb{L}^2(\mu)$,
\begin{equation}
\label{2.5}
\lim_{t\rightarrow\infty}e^{\lambda_1t}\langle g,P_tf\rangle_\mu=\langle\eta_1,f\rangle_\mu\langle\eta_1,g\rangle_\mu.
\end{equation}
\par $(\mathrm{ii})$ for any $n\geq1, \eta_n\in \mathbb{L}^1(\mu)$.
\par $(\mathrm{iii})$ for all $x > 0$ and all $t > 0$, there exists some
density $r(t,x,\cdot)$ that satisfies
\begin{equation}
\label{2.6}
\mathbb{E}_x[f(X_t), T_{0}>t]=\int_{0}^{\infty}r(t,x,y)f(y)\mu(dy)
\end{equation}
for all bounded Borel function $f$. Moreover, the density
$r(t,x,\cdot)\in\mathbb{L}^2(\mu)$ for all $x > 0, t > 0$. In particular, there exists a function $B(t)\geq0, \lim_{t\rightarrow\infty}B(t)=0$ such that
\begin{equation}
\label{2.7}
\int_{0}^{\infty}r^2(t,x,y)\mu(dy)<r(t,x,x)B(t)<\infty.
\end{equation}
\par $(\mathrm{iv})$ for all $x > 0$ and $A\in\mathcal{B}(0,\infty)$,
\begin{equation}
\label{2.8}
\lim_{t\rightarrow\infty}e^{\lambda_1t}\mathbb{P}_{x}(T_{0}>t)=\eta_{1}(x)\langle\eta_{1},1\rangle_\mu.
\end{equation}
\end{prop}
\par
The following result plays an important role in our following arguments and is of independent interest.
\begin{prop}
\label{prop 2.3}
Assume $(\mathrm{H})$ holds. Then $\eta_1$ is bounded.
\end{prop}
\begin{pf}
First, if $(\mathrm H)$ is satisfied, based on \cite[Proposition 7.6]{CCLMMS09}, we can deduce that there is $x_0>0$ such that $B_1:=\sup_{x\geq x_0}\mathbb{E}_x[e^{\lambda_1T_{x_0}}]<\infty$. From the equality (\ref{2.8}), we get that $B_2:=\sup_{u\geq 0}e^{\lambda_1u}\mathbb{P}_{x_0}(T_{0}>u)<\infty$. Then for $x>x_0$, we have
\begin{eqnarray*}
\mathbb{P}_{x}(T_{0}>t)&=&\int_{0}^{t}\mathbb{P}_{x_0}(T_{0}>u)\mathbb{P}_{x}(T_{x_0}\in d(t-u))+\mathbb{P}_{x}(T_{x_0}>t)\\
&\leq&B_2\int_{0}^{t}e^{-\lambda_1u}\mathbb{P}_{x}(T_{x_0}\in d(t-u))+\mathbb{P}_{x}(T_{x_0}>t)\\
&\leq&B_2e^{-\lambda_1t}\mathbb{E}_x[e^{\lambda_1T_{x_0}}]+e^{-\lambda_1t}\mathbb{E}_x[e^{\lambda_1T_{x_0}}]\\
&\leq&e^{-\lambda_1t}B_1(B_2+1).
\end{eqnarray*}
Thus, we get that $e^{\lambda_1t}\mathbb{P}_x(T_{0}>t)$ is uniformly bounded in the variables $t$ and $x$. By using the equality (\ref{2.8}) again, it is easily seen that for $x>x_0>0$, $\eta_1(x)\leq\frac{B_1(B_2+1)}{\langle\eta_1,1\rangle_\mu}$. On the other hand, for $0<x\leq x_0$, we have $\mathbb{P}_{x}(T_{0}>t)\leq\mathbb{P}_{x_0}(T_{0}>t)$. Thus from the equality (\ref{2.8}), we get that $\eta_1(x)\leq\eta_1(x_0)$. Hence, for any $x>0$, there exists $x_0>0$ such that $\eta_1(x)\leq\max\{\eta_1(x_0),\frac{B_1(B_2+1)}{\langle\eta_1,1\rangle_\mu}\}$. This completes the proof.
\qed
\end{pf}

\section{Quasi-ergodic distribution}
\label{sect3}
In this section, we study the existence and uniqueness of quasi-ergodic distribution for one-dimensional diffusion $X$ killed at 0, when 0 is an exit boundary and $+\infty$ is an entrance boundary. More formally, the following definition captures the main object of interest of this work.
\begin{defi}
\label{defi3.1}
We say that $\nu \in\mathscr{P}(0,\infty)$ is a quasi-ergodic distribution if, for all $t>0$ and any $A\in\mathcal{B}(0,\infty)$, there exists a $\pi\in\mathscr{P}(0,\infty)$ such that the following limit exists in the weak sense
\begin{equation*}
\lim_{t\rightarrow\infty}\mathbb{E}_\pi(\frac{1}{t}\int_{0}^{t}{\bf1}_{A}(X_s)ds|T_{0}>t)=\nu(A).
\end{equation*}
\end{defi}
\par
Let us recall the notion of intrinsic ultracontractivity $\mathrm{[IU]}$, which was introduced by Davies and Simon \cite{DS84}, is a very important concept in both analysis and probability and has been studied extensively. In addition we point out that for one-dimensional generalized diffusion processes with no natural boundary, a sufficient condition for $\mathrm{[IU]}$ was given in \cite[Theorem 2.11]{T07}
\begin{defi}
The semigroup $\left(P_t\right)_{t\geq0}$ is said to be intrinsically ultracontractive $\mathrm{[IU]}$ if, for any $t>0$, there exists a constant $c_t>0$ such that
\begin{equation}
\label{3.1}
r(t,x,y)\leq c_t\eta_1(x)\eta_1(y)~~~~\mathrm{for}~x,~y\in (0,\infty).
\end{equation}
\end{defi}
\par
Define
\begin{equation*}
\nu_1(A):=\int_A\eta^2_1(x)\mu(dx),~~~~~~~~~~A\in\mathcal{B}(0,\infty).
\end{equation*}
We know from Proposition \ref{prop 2.2} that $\|\eta_1\|_{\mathbb{L}^2(\mu)}=1$, then $\nu_1$ is a distribution on $(0,\infty)$.
\par We may now state the following result.
\begin{thm}
\label{thm3.1}
Assume that $(\mathrm{H})$ and $\mathrm{[IU]}$ are satisfied. Then for any $0<q<1$, any $\nu \in\mathscr{P}(0,\infty)$ and any bounded Borel function $f$ on $(0,\infty)$, we have
\begin{equation*}
\lim_{t\rightarrow\infty}\mathbb{E}_{\nu}(f(X_{qt})|T_{0}>t)=\int_{0}^{\infty}f(y)\nu_1(dy).
\end{equation*}
\end{thm}
\begin{pf}
We know from Proposition \ref{prop 2.2} that $r(1,x,\cdot)\in\mathbb{L}^2(\mu)$. By using Proposition \ref{prop 2.2}, writing $r(t,x,\cdot)=P_{t-1}r(1,x,\cdot)$ $\mu$-a.s. and noticing that
\begin{equation*}
\int_{0}^{\infty} r(1,x,y)\eta_1(y)\mu(dy)=(P_1\eta_1)(x)=e^{-\lambda_1}\eta_1(x),
\end{equation*}
we deduce that the following limit exists in $\mathbb{L}^2(\mu)$
\begin{equation}
\label{3.2}
\lim\limits_{t\rightarrow\infty}e^{\lambda_1t}r(t,x,\cdot)=e^{\lambda_1}\langle r(1,x,\cdot),\eta_1\rangle_\mu\eta_1(\cdot)=\eta_1(x)\eta_1(\cdot).
\end{equation}
Similarly, for $t>\min\{\frac{1}{q}, \frac{1}{1-q}\}$, writing $r(qt,x,\cdot)=P_{qt-1}r(1,x,\cdot)$ $\mu$-a.s. and
$r(t-qt,z,\ast)=P_{t-qt-1}r(1,z,\ast)$ $\mu$-a.s., we deduce that the following limit exists in $\mathbb{L}^2(\mu)$
\begin{equation}
\label{3.3}
\lim\limits_{t\rightarrow\infty}e^{\lambda_1t}r(qt,x,\cdot)r(t-qt,z,\ast)=e^{2\lambda_1}\langle r(1,x,\cdot),\eta_1\rangle_\mu\eta_1(\cdot)\langle r(1,z,\ast),\eta_1\rangle_\mu\eta_1(\ast)=\eta_1(x)\eta_1(\cdot)\eta_1(z)\eta_1(\ast).
\end{equation}
By the semigroup property, for any $(t,x,y)\in(0,\infty)\times(0,\infty)\times(0,\infty)$, we have
\begin{equation*}
r(t,x,y)=\int_{0}^{\infty}r(qt,x,z)r(t-qt,z,y)\mu(dz).
\end{equation*}
We know from the definition of [IU] that the equality $(\ref{3.1})$ holds, then for $t>1, e^{\lambda_1t}r(t,x,\cdot)\in \mathbb{L}^1(\mu)$ and is dominated by $e^{\lambda_1t}c_t\eta_1(x)\eta_1(\cdot)$. Similarly, the equality $(\ref{3.1})$ holds, then for $t>\min\{\frac{1}{q}, \frac{1}{1-q}\}$, $e^{\lambda_1t}r(qt,x,\cdot)r(t-qt,z,\ast)$ is dominated by $e^{\lambda_1t}c_{qt}c_{t-qt}\eta_1(x)\eta_1(z)\eta_1(\cdot)\eta_1(\ast)$. Note that $\eta_{1}$ is bounded (see Proposition \ref{prop 2.3}) and $\eta_{1}\in \mathbb{L}^1(\mu)$
(see Proposition \ref{prop 2.2}). It now follows from the Lebesgue dominated convergence theorem and Fubini's theorem that for any $0<q<1$ and any bounded Borel function $f$ on $(0,\infty)$,
\begin{eqnarray*}
\lim_{t\rightarrow\infty}\mathbb{E}_{\nu}(f(X_{qt})|T_{0}>t)
&=&\lim_{t\rightarrow\infty}\frac{\int_{0}^{\infty}\mathbb{E}_{x}(f(X_{qt}),T_{0}>t)\nu(dx)}{\int_{0}^{\infty}\mathbb{P}_{x}(T_{0}>t)\nu(dx)}\\
&=&\lim_{t\rightarrow\infty}\frac{\int_{0}^{\infty}\int_{0}^{\infty}\int_{0}^{\infty}r(qt,x,z)f(z)r(t-qt,z,y)\mu(dz)\mu(dy)\nu(dx)}{\int_{0}^{\infty}\int_{0}^{\infty}r(t,x,y)\mu(dy)\nu(dx)}\\
&=&\lim_{t\rightarrow\infty}\frac{e^{\lambda_1t}\int_{0}^{\infty}\int_{0}^{\infty}\int_{0}^{\infty}r(qt,x,z)f(z)r(t-qt,z,y)\mu(dz)\mu(dy)\nu(dx)}{e^{\lambda_1t}\int_{0}^{\infty}\int_{0}^{\infty}r(t,x,y)\mu(dy)\nu(dx)}\\
&=&\frac{\int_{0}^{\infty}\int_{0}^{\infty}\int_{0}^{\infty}\eta_1(x)\eta^2_1(z)\eta_1(y)f(z)\mu(dz)\mu(dy)\nu(dx)}{\int_{0}^{\infty}\int_{0}^{\infty}\eta_1(x)\eta_1(y)\mu(dy)\nu(dx)}\\
&=&\int_{0}^{\infty}\eta^2_1(z)f(z)\mu(dz)\\
&=&\int_{0}^{\infty}f(y)\nu_1(dy).
\end{eqnarray*}
\qed
\end{pf}
\par
Notice that when $q=1$, we know from \cite[Theorem 4.2]{LJ12} that for $A\in\mathcal{B}(0,\infty), \lim_{t\rightarrow\infty}\mathbb{P}_{x}(X_t\in A|T_{0}>t)=\nu_2(A)$, where $\nu_2(A)=\frac{\int_A\eta_1(x)\mu(dx)}{\int_{0}^{\infty}\eta_1(x)\mu(dx)}$ is a quasi-stationary distribution. As we have seen that the above result exhibits a phase transition.
A further investigation of $\nu_1$ shows that $\nu_1$ can also be described as the following double limit.
\begin{thm}
\label{thm3.2}
Assume that $(\mathrm{H})$ and $\mathrm{[IU]}$ are satisfied. Then for any $\nu \in\mathscr{P}(0,\infty)$ and any bounded Borel function $f$ on $(0,\infty)$, we have
\begin{eqnarray*}
\lim_{t\rightarrow\infty}\lim_{T\rightarrow\infty}\mathbb{E}_{\nu}(f(X_{t})|T_{0}>T)=\int_{0}^{\infty}f(y)\nu_1(dy).
\end{eqnarray*}
\end{thm}
\begin{pf}
For $t>1, T-t>1$, writing $r(T-t,z,\cdot)=P_{T-t-1}r(1,z,\cdot)$ $\mu$-a.s., we deduce that the following limit exists in $\mathbb{L}^2(\mu)$
\begin{equation}
\label{3.4}
\lim\limits_{T\rightarrow\infty}e^{\lambda_1T}r(T-t,z,\cdot)=e^{\lambda_1(t+1)}\langle r(1,z,\cdot),\eta_1\rangle_\mu\eta_1(\cdot)=e^{\lambda_1t}\eta_1(z)\eta_1(\cdot).
\end{equation}
Moreover, we know from the definition of [IU] that the equality $(\ref{3.1})$ holds, then for $t>1, T-t>1$, $e^{\lambda_1t}r(t,x,\cdot)$ is dominated by $e^{\lambda_1t}c_t\eta_1(x)\eta_1(\cdot)$, $e^{\lambda_1T}r(T,x,\cdot)$ is dominated by $e^{\lambda_1T}c_T\eta_1(x)\eta_1(\cdot)$ and $e^{\lambda_1T}r(t,x,\cdot)r(T-t,z,\ast)$ is dominated by $e^{\lambda_1T}c_{t}c_{T-t}\eta_1(x)\eta_1(z)\eta_1(\cdot)\eta_1(\ast)$.
Then by the Lebesgue dominated convergence theorem, the semigroup property, the equalities (\ref{3.2}), (\ref{3.4}) and Fubini's theorem, we have
\begin{eqnarray*}
\lim_{t\rightarrow\infty}\lim_{T\rightarrow\infty}\mathbb{E}_{\nu}(f(X_{t})|T_{0}>T)&=&\lim_{t\rightarrow\infty}\lim_{T\rightarrow\infty}\frac{\int_{0}^{\infty}\mathbb{E}_{x}(f(X_{t}),T_{0}>T)\nu(dx)}{\int_{0}^{\infty}\mathbb{P}_{x}(T_{0}>T)\nu(dx)}\\
                                                                                              &=&\lim_{t\rightarrow\infty}\lim_{T\rightarrow\infty}\frac{\int_{0}^{\infty}\int_{0}^{\infty}\int_{0}^{\infty}r(t,x,z)f(z)r(T-t,z,y)\mu(dz)\mu(dy)\nu(dx)}{\int_{0}^{\infty}\int_{0}^{\infty}r(T,x,y)\mu(dy)\nu(dx)}\\
                                                                                              &=&\lim_{t\rightarrow\infty}\lim_{T\rightarrow\infty}\frac{e^{\lambda_1T}\int_{0}^{\infty}\int_{0}^{\infty}\int_{0}^{\infty}r(t,x,z)f(z)r(T-t,z,y)\mu(dz)\mu(dy)\nu(dx)}{e^{\lambda_1T}\int_{0}^{\infty}\int_{0}^{\infty}r(T,x,y)\mu(dy)\nu(dx)}\\
                                                                                              &=&\lim_{t\rightarrow\infty}\frac{\int_{0}^{\infty}\int_{0}^{\infty}\int_{0}^{\infty}e^{\lambda_1t}r(t,x,z)f(z)\eta_1(z)\eta_1(y)\mu(dz)\mu(dy)\nu(dx)}{\int_{0}^{\infty}\int_{0}^{\infty}\eta_1(x)\eta_1(y)\mu(dy)\nu(dx)}\\
                                                                                              &=&\lim_{t\rightarrow\infty}\frac{\int_{0}^{\infty}\int_{0}^{\infty}e^{\lambda_1t}r(t,x,z)f(z)\eta_1(z)\mu(dz)\nu(dx)}{\int_{0}^{\infty}\eta_1(x)\nu(dx)}\\
                                                                                              &=&\frac{\int_{0}^{\infty}\int_{0}^{\infty}\eta_1(x)\eta^2_1(z)f(z)\mu(dz)\nu(dx)}{\int_{0}^{\infty}\eta_1(x)\nu(dx)}\\
                                                                                              &=&\int_{0}^{\infty}\eta^2_1(z)f(z)\mu(dz)=\int_{0}^{\infty}f(y)\nu_1(dy).
\end{eqnarray*}
\qed
\end{pf}
\par
We now describe the quasi-ergodic behavior of one-dimensional diffusion $X$ killed at 0, when 0 is an exit boundary and $+\infty$ is an entrance boundary. The following result implies that $\nu_1$ is the unique quasi-ergodic distribution of the process $X$.
\begin{thm}
\label{thm3.3}
Assume that $(\mathrm{H})$ and $\mathrm{[IU]}$ are satisfied. Then for any $\nu \in\mathscr{P}(0,\infty)$ and any bounded Borel function $f$ on $(0,\infty)$, we have
\begin{eqnarray*}
\lim_{t\rightarrow\infty}\mathbb{E}_{\nu}(\frac{1}{t}\int_{0}^{t}f(X_{s})ds|T_{0}>t)=\int_{0}^{\infty}f(y)\nu_1(dy).
\end{eqnarray*}
In particular, $\nu_1$ is just the stationary distribution of $Q$-process$($the process conditioned to never be absorbed$)$.
\end{thm}
\begin{pf}
We set $s=qt$. Then by change of variable in the Lebesgue integral, the Lebesgue dominated convergence theorem and Theorem \ref{thm3.1}, we obtain
\begin{eqnarray*}
\lim_{t\rightarrow\infty}\mathbb{E}_{\nu}(\frac{1}{t}\int_{0}^{t}f(X_{s})ds|T_{0}>t)
&=&\lim_{t\rightarrow\infty}\mathbb{E}_{\nu}(\int_{0}^{1}f(X_{qt})dq|T_{0}>t)\\
&=&\lim_{t\rightarrow\infty}\int_{0}^{1}\mathbb{E}_{\nu}(f(X_{qt})|T_{0}>t)dq\\
&=&\int_{0}^{\infty}f(y)\nu_1(dy).
\end{eqnarray*}
For the last part of the theorem, it can be found in \cite[Corollary 6.2]{CCLMMS09}.
\qed
\end{pf}
\par
Next, we introduce the ultracontractivity of the semigroup, which is of independent interest. We say that a semigroup ${P_t}$ is called ultracontractive if $P_t$ are bounded operators from $\mathbb{L}^1$ to $\mathbb{L}^{\infty}$ for all $t>0$. If a semigroup ${P_t}$ is ultracontractive, it can deduce that for $t>0$, there exists positive constant $c_t$ such that
\begin{equation*}
\label{3.3}
r(t,x,y)\leq c_t<\infty~~~~\mathrm{for}~x,~y\in (0,\infty).
\end{equation*}
\begin{prop}
Assuming that absorption is certain for the process $X$, i.e. $(\mathrm{H1})$ holds. Then the following are equivalent$:$
\par $(\mathrm{i})$ $(P_t)_{t\geq0}$ is ultracontractive$;$
\par $(\mathrm{ii})$ for $q>2$, the following conditions hold$:$
\begin{equation}
\label{3.5}
\sup\limits_{0<x<1}\left(\Lambda(x)\right)^{q/(q-2)}\mu([x,1))<+\infty,
\end{equation}
\begin{equation}
\label{3.6}
\sup\limits_{x\geq1}\left(\Lambda(x)\right)^{q/(q-2)}\mu([x,+\infty))<+\infty.
\end{equation}
\end{prop}
\begin{pf}
We consider $Z_{t}=\Lambda(X_t)$. It is direct to prove that $Z$ is in natural scale on the interval $(\Lambda(0),\infty)$, that is, for
$\Lambda(0)<a\leq y\leq b<\infty=\Lambda(\infty)$
\begin{equation*}
\mathbb{P}_y\left(T^{Z}_{a}<T^{Z}_{b}\right)=\frac{b-y}{b-a},
\end{equation*}
where $T^{Z}_{a}$ be the hitting time of $a$ for the process $Z$. We know from \cite[Theorem 3.5]{S14} that $(P^{Z}_t)_{t\geq0}$ is ultracontractive for $Z$ if and only if, for $q>2$ the following conditions hold:
\begin{equation}
\label{3.7}
\sup\limits_{0<z<1}z^{q/(q-2)}m([z,1))<+\infty,
\end{equation}
\begin{equation}
\label{3.8}
\sup\limits_{z\geq1}z^{q/(q-2)}m([z,+\infty))<+\infty,
\end{equation}
where $m$ is the speed measure of $Z$, which is given by
\begin{equation*}
m(dz)=\frac{2dz}{(\Lambda'(\Lambda^{-1}(z)))^{2}}
\end{equation*}
(see \cite[formula (5.51)]{KS88}), because $Z$ satisfies the SDE
\begin{equation*}
dZ_{t}=\Lambda'(\Lambda^{-1}(Z_{t}))dB_{t}.
\end{equation*}
After a change of variables, we obtain
\begin{equation*}
\label{}
\sup\limits_{0<z<1}z^{q/(q-2)}m([z,1))=2\sup\limits_{0<x<1}\left(\Lambda(x)\right)^{q/(q-2)}\mu([x,1)),
\end{equation*}
\begin{equation*}
\label{}
\sup\limits_{z\geq1}z^{q/(q-2)}m([z,+\infty))=2\sup\limits_{x\geq1}\left(\Lambda(x)\right)^{q/(q-2)}\mu([x,+\infty)).
\end{equation*}
Therefore, we just prove the equivalence.
\qed
\end{pf}

\section{An example}
\label{sect4}
We consider a famous biological model: the logistic Feller diffusion process
\begin{equation*}
dY_t=\sqrt{\gamma Y_{t}}dB_t+\left(rY_{t}-cY^{2}_{t}\right)dt,~~~~~~~~~~Y_{0}=y>0,
\end{equation*}
where $\gamma$, $r$ and $c$ are assumed to be positive constants. After a suitable change of variables: $X_t=2\sqrt{Y_{t}/\gamma}$, we obtain
\begin{equation*}
\alpha(x)=\frac{1}{2x}-\frac{rx}{2}+\frac{c\gamma x^{3}}{8}.
\end{equation*}
According to the proof of \cite[Theorem 8.2]{CCLMMS09} and \cite[Theorem 3.4]{M14}, we know that the hypotheses $(\mathrm{H})$ and $\mathrm{[IU]}$ are satisfied for $X$ respectively.
So it follows from the previous section that there exists a unique quasi-ergodic distribution for $Y$.

\section*{Acknowledgements}
The work is supported by Hunan Provincial Innovation Foundation For Postgraduate (Grant No. CX2015B203) and the National Natural Science Foundation of China (Grant No.11371301).


\begin{thebibliography}{00}

\bibitem{BR99}\label{BR99}
Breyer, L. A., and Roberts, G. O. 1999.
A quasi-ergodic theorem for evanescent processes.
{\em Stochastic Processes. Appl.} 84: 177--186.

\bibitem{CCLMMS09}\label{CCLMMS09}
Cattiaux, P., Collet, P., Lambert, A., Mart\'{\i}nez, S., M\'{e}l\'{e}ard, S., and San Mart\'{\i}n, J. 2009.
Quasi-stationary distributions and diffusion models in population dynamics.
{\em Ann. Probab.} 37: 1926--1969.

\bibitem{CJ13}\label{CJ13}
Chen, J. W., and Jian, S. Q. 2013.
Some limit theorems of killed Brownian motion.
{\em Sci. China. Math.} 56: 497--514.

\bibitem{CJ14}\label{CJ14}
Chen, J. W., and Jian, S. Q. 2014.
A remark on quasi-ergodicity of ultracontractive Markov processes.
{\em Statist. Probab. Lett.} 87: 184--190.

\bibitem{CLJ12}\label{CLJ12}
Chen, J. W., Li, H. T., and Jian, S. Q. 2012.
Some limit theorems for absorbing Markov processes.
{\em J. Phys. A: Math. Theor.} 45: 345003 (11pp).

\bibitem{DS84}\label{DS84}
Davies, E. B., and Simon, B. 1984.
Ultracontractivity and the heat kernel for Schr\"{o}inger operators and Dirichlet Laplacians.
{\em J. Funct. Anal.} 59: 335--395.

\bibitem{IW89}\label{IW89}
Ikeda, N., and Watanabe, S. 1989.
{\em Stochastic Differential Equations and Diffusion Processes} (North-Holland Mathematical Library {\bf24}), 2nd edn. North-Holland, Amsterdam.

\bibitem{KS88}\label{KS88}
Karatzas, I., and Shreve, S. 1988.
{\em Brownian Motion and Stochastic Calculus.} Springer, New York.

\bibitem{KT81}\label{KT81}
Karlin, S., and Taylor, H. M. 1981.
{\em A Second Course in Stochastic Processes,} 2nd edn. Academic Press, New York.

\bibitem{LJ12}\label{LJ12}
Littin, J. 2012.
Uniqueness of quasistationary distributions and discrete spectra when $\infty$ is an entrance boundary and 0 is singular.
{\em J. Appl. Probab.} 49: 719--730.

\bibitem{M14}\label{M14}
Miura, Y. 2014.
Ultracontractivity for Markov semigroups and quasi-stationary distributions.
{\em Stoch. Anal. Appl.} 32: 591--601.

\bibitem{S14}\label{S14}
Shigekawa, I. 2014.
Dual ultracontractivity and its applications.
{\em Front. Math. China.} 9: 899--928.

\bibitem{T07}\label{T07}
Tomisaki, M. 2007.
Intrinsic ultracontractivity and small perturbation for one-dimensional generalized diffusion operators.
{\em J. Funct. Anal.} 251: 289--324.

\bibitem{ZLS14}\label{ZLS14}
Zhang, J. F., Li, S. M., and Song, R. M. 2014.
Quasi-stationarity and quasi-ergodicity of general Markov processes.
{\em Sci. China. Math.} 57: 2013--2024.


\end{thebibliography}
\end{document}